\documentclass[11pt]{article}
\usepackage{amsthm,amstext}
\usepackage{amsmath}
\usepackage{graphicx}
\usepackage{amsfonts}
\usepackage{amssymb}
\theoremstyle{plain}
\newtheorem{theorem}{Theorem}[section]

\newtheorem{proposition}[theorem]{Proposition}

\theoremstyle{definition}
\newtheorem{definition}[theorem]{Definition}
\newtheorem{corollary}[theorem]{Corollary}
\newtheorem{note}[theorem]{Note}
\newtheorem{notation}[theorem]{Notation}
\newtheorem{example}{\sc Example}
\theoremstyle{remark}
\newtheorem{remark}{\sc Remark}

\date{}
\title{\bf On Intuitionistic Fuzzy Magnified Translation \\in Semigroups}\vspace{.25 in}
\author{ {\bf Sujit Kumar Sardar$^1$}\\\\\\ {\bf Manasi Mandal$^2$} \\\\and \\\\{\bf Samit Kumar Majumder$^3$}\\
Department of Mathematics, Jadavpur\\
University, Kolkata-700032, INDIA\\
{\tt $^1$sksardarjumath@gmail.com}\\
{\tt $^2$}manasi$_{-}$ju@yahoo.in\\
{\tt $^3$samitmathalg@yahoo.co.in}
 }

\begin{document}
\maketitle

\begin{abstract}

The notion of intuitionistic fuzzy sets was introduced by {\it Atanassov} as a generalization of the notion of fuzzy sets. {\it S.K Sardar} and {\it S.K. Majumder} unified the idea of fuzzy translation and fuzzy multiplication of {\it Vasantha Kandasamy} to introduce the concept of fuzzy magnified translation in groups and semigroups. The purpose of this paper is to intuitionistically fuzzify$($by using {\it Atanassov's} idea$)$ the concept of fuzzy magnified translation in semigroups. Here among other results we obtain some characterization theorems of regular, intra-regular, left$($right$)$ regular semigroups in terms of intuitionistic fuzzy magnified translation.\\

\end{abstract}

\section{Introduction}
A  semigroup is an algebraic structure consisting of a non-empty set $S$ together with an associative binary operation\cite{H}. The formal study of semigroups began in the early 20th century. Semigroups are important in many areas of mathematics, for example, coding and language theory, automata theory, combinatorics and mathematical analysis. The concept of fuzzy sets was introduced by {\it Lofti Zadeh}\cite{Z} in his classic paper in $1965.$ {\it Azirel Rosenfeld}\cite{R} used the idea of fuzzy set to introduce the notions of fuzzy subgroups. {\it Nobuaki Kuroki}\cite{K1,K2,K3} is the pioneer of fuzzy ideal theory of semigroups. The idea of fuzzy subsemigroup was also introduced by {\it Kuroki}$\cite{K1,K3}$. In $\cite{K2}$, {\it Kuroki} characterized several classes of semigroups in terms of fuzzy left, fuzzy right and fuzzy bi-ideals. Others who worked on fuzzy semigroup theory, such as {\it X.Y. Xie}\cite{X1,X2}, {\it Y.B. Jun}\cite{J1,J2}, are mentioned in the bibliography. The notion of intuitionistic fuzzy sets was introduced by {\it Atanassov}\cite{A1,A2} as a generalization of the notion of fuzzy sets. In this paper we introduce the notion of intuitionistic fuzzy magnified translation in semigroups and observe some of its important properties. Here we characterize regular, intra-regular and left$($right$)$ regular semigroups in terms of intuitionistic fuzzy magnified translation. Finally we also observe that intuitionistic fuzzy translation and intuitionistic fuzzy multiplication are the particular cases of intuitionistic fuzzy magnified translation.
\section{Preliminaries}

\indent In this section we discuss some elementary definitions that we use in the sequel.\\

If $(X,\ast)$ is a mathematical system such that $\forall a,b,c\in X,$ $(a\ast b)\ast c=a\ast(b\ast c),$ then $\ast$ is called associative and $(X,\ast)$ is called a {\it semigroup}\cite{M}.\\

Throughout the paper unless otherwise stated $S$ will denote a semigroup.\\

A non-empty subset $A$ of a semigroup $S$ is called a {\it subsemigroup}\cite{K2} of $S$ if $AA\subseteq A.$

A subsemigroup $A$ of $S$ is called a {\it bi-ideal}\cite{K2} of $S$ if $ASA\subseteq A.$

A subsemigroup $A$ of $S$ is called an {\it $(1,2)$-ideal}\cite{K2} of $S$ if $ASAA\subseteq A.$

A {\it left} ({\it right})
{\it ideal}\cite{M} of a semigroup $S$ is a non-empty subset $I$ of $S$ such that $SI \subseteq I$ ($IS \subseteq I$). If $I$ is both a left and a right ideal of a semigroup $S$, then we say that $I$ is an {\it ideal}\cite{M} of $S$.\\

A {\it fuzzy subset}\cite{Z} in $S$ is a function $\mu : S \longrightarrow [0,1]$.

\indent Let $\mu $ be a fuzzy subset of a set $X$ and $\alpha \in\lbrack 0,1-\sup \{\mu (x):x\in X\}].$ A mapping $\mu_{\alpha }^{T}:X\rightarrow \lbrack 0,1]$\ is called a {\it fuzzy translation}\cite{K} of $\mu $
if $\mu _{\alpha }^{T}(x)=\mu (x)+\alpha $ for all $x\in X$.

\indent Let $\mu $ be a fuzzy subset of a set $X$ and $\beta \in\lbrack 0,1].$ A mapping $\mu _{\beta }^{M}:X\rightarrow \lbrack 0,1]$\ is called a {\it fuzzy multiplication}\cite{K} of $\mu $ if $\mu _{\beta }^{M}(x)=\beta\cdot\mu (x)$ for all $x\in X$.

Let $\mu $ be a fuzzy subset of a set $X,$ $\alpha \in\lbrack 0,1-\sup \{\mu (x):x\in X\}]$ and $\beta\in \lbrack 0,1].$ A mapping $\mu _{\beta\alpha }^{C}:X\rightarrow \lbrack 0,1]$\ is called a {\it fuzzy magnified translation}\cite{S1} of $\mu $ if $\mu _{\beta\alpha }^{C}(x)=\beta\cdot\mu (x)+\alpha $ for
all $x\in X$.

A non-empty fuzzy subset $\mu$ of a semigroup $S$ is called a {\it fuzzy subsemigroup}\cite{M} of $S$ if $\mu(xy)\geq\min\{\mu(x),\mu(y)\}\forall x,y\in S.$

A fuzzy subsemigroup $\mu$ of a semigroup $S$ is called a {\it fuzzy bi-ideal}\cite{M} of $S$ if $\mu(xyz)\geq\min\{\mu(x),\mu(z)\}\forall x,y,z\in S.$

A fuzzy subsemigroup $\mu$ of a semigroup $S$ is called a {\it fuzzy $(1,2)$-ideal}\cite{M} of $S$ if $\mu(x\omega(yz))\geqslant\min\{\mu(x),\mu(y),\mu(z)\}\forall x,\omega,y,z\in S.$

A non-empty fuzzy subset $\mu$ of a semigroup $S$ is called a {\it fuzzy left$($right$)$ ideal}\cite{M} of $S$ if $\mu(xy)\geq\mu(y)($resp. $\mu(xy)\geq\mu(x))$\ $\forall x,y\in S.$

A non-empty fuzzy subset $\mu$ of a semigroup $S$ is called a {\it fuzzy two-sided ideal} or a {\it fuzzy ideal}\cite{M} of $S$ if it is both a fuzzy left and a fuzzy right ideal of $S.$

A fuzzy ideal $\mu$ of a semigroup $S$ is called a {\it fuzzy semiprime ideal}\cite{K1} of $S$ if $\mu(x)\geq\mu(x^{2})\forall x\in S.$\\

{\it Atanassov} introduced in \cite{A1,A2} the concept of intuitionistic fuzzy sets defined on a non-empty set $X$ as objects having the form\\
$$A=\{<x,\mu_{A}(x),\nu_{A}(x)>:x\in X\},$$ where the functions $\mu_{A}: X\rightarrow [0,1]$ and $\nu_{A}: X\rightarrow [0,1]$ denote the degree of membership and the degree of non-membership of each element $x\in X$ to the set $A$ respectively, and $0\leq \mu_{A}(x)+\nu_{A}(x)\leq 1$ for all $x\in X.$\\
\indent Such defined objects are studied by many authors and have many interesting applications in mathematics.\\
\indent Let $A$ and $B$ be two intuitionistic fuzzy subsets of a set $X.$ Then the following expressions are defined in \cite{A1,A2}.\\
\indent $(i)$ $A\subseteq B$ if and only if $\mu_{A}(x)\leq\mu_{B}(x)$ and $\nu_{A}(x)\geq\nu_{B}(x),$\\
\indent $(ii)$ $A=B$ if and only if $A\subseteq B$ and $B\subseteq A,$\\
\indent $(iii)$ $A^{C}=\{<x,\nu_{A}(x),\mu_{A}(x)>:x\in X\},$\\
\indent $(iv)$ $A\cap B=\{<x,\min\{\mu_{A}(x),\mu_{B}(x)\},\max\{\nu_{A}(x),\nu_{B}(x)\}>:x\in X\},$\\
\indent $(v)$ $A\cup B=\{<x,\max\{\mu_{A}(x),\mu_{B}(x)\},\min\{\nu_{A}(x),\nu_{B}(x)\}>:x\in X\}.$\\

From the definition it follows that $A\cap B$ is the same as $\mu_{A}\cap\mu_{B}$ and $\nu_{A}\cup\nu_{B}.$ Also $A\cup B$ is the same as $\mu_{A}\cup\mu_{B}$ and $\nu_{A}\cap\nu_{B}.$\\

For the sake of simplicity, we shall use the symbol $A=(\mu_{A},\nu_{A})$ for the intuitionistic fuzzy subset $A=\{<x,\mu_{A}(x),\nu_{A}(x)>:x\in X\}.$\\

A non-empty intuitionistic fuzzy subset $A=(\mu_{A},\nu_{A})$ of a semigroup $S$ is
called an {\it intuitionistic fuzzy subsemigroup} of $S$ if $(i)$ $\mu_{A}(xy)\geq\min\{\mu_{A}(x),\mu_{A}(y)\}\forall x,y\in S,$ $(ii)$ $\nu_{A}(xy)\leq\max\{\nu_{A}(x),\nu_{A}(y)\}\forall x,y\in S.$

An intuitionistic fuzzy subsemigroup $A=(\mu_{A},\nu_{A})$ of a semigroup $S$ is called an {\it intuitionistic fuzzy bi-ideal} of $S$ if $(i)$ $\mu_{A}(xyz)\geq\min\{\mu_{A}(x),\mu_{A}(z)\}\forall x,y,z\in S,$ $(ii)$ $\nu_{A}(xyz)\leq\max\{\nu_{A}(x),\nu_{A}(z)\}\forall x,y,z\in S.$

An intuitionistic fuzzy subsemigroup $A=(\mu_{A},\nu_{A})$ of a semigroup $S$ is called an {\it intuitionistic fuzzy $(1,2)$-ideal} of $S$ if $(i)$ $\mu_{A}(x\omega(yz))\geqslant\min\{\mu_{A}(x),\mu_{A}(y),\mu_{A}(z)\}\forall
x,\omega,y\\,z\in S,$ $(ii)$ $\nu_{A}(x\omega(yz))\leqslant\max\{\nu_{A}(x),\nu_{A}(y),\nu_{A}(z)\}\forall
x,\omega,y,z\in S.$

A non-empty intuitionistic fuzzy subset $A=(\mu_{A},\nu_{A})$ of a semigroup $S$ is called an {\it intuitionistic fuzzy left$($right$)$ ideal} of $S$ if $(i)$ $\mu_{A}(xy)\geq\mu(y)($resp. $\mu_{A}(xy)\geq\mu_{A}(x))$\ $\forall x,y\in S,$ $(ii)$ $\nu_{A}(xy)\leq\nu_{A}(y)($resp. $\nu_{A}(xy)\leq\nu_{A}(x))$\ $\forall x,y\in S.$

A non-empty intuitionistic fuzzy subset $A=(\mu_{A},\nu_{A})$ of a semigroup $S$ is called an {\it intuitionistic fuzzy two-sided ideal} or an {\it intuitionistic fuzzy ideal} of $S$ if it is both an intuitionistic fuzzy left and an intuitionistic fuzzy right ideal of $S.$

An intuitionistic fuzzy ideal $A=(\mu_{A},\nu_{A})$ of a semigroup $S$ is called an {\it intuitionistic fuzzy semiprime ideal} of $S$ if $(i)$ $\mu_{A}(x)\geq\mu_{A}(x^{2})\forall x\in S,$ $(ii)$ $\nu_{A}(x)\leq\nu_{A}(x^{2})\forall x\in S.$

\section{Main Results}

\begin{definition}
Let $A=(\mu_{A},\nu_{A}) $ be an intuitionistic fuzzy subset of a set $X,$ $\alpha \in\lbrack 0,\inf \{\nu_{A}(x):x\in X\}].$ An object having the form $A^{T}_{\alpha}=((\mu_{A})^{T}_{\alpha},(\nu_{A})^{T}_{\alpha})$ is called an {\it intuitionistic fuzzy translation} of $A$ if $(\mu_{A})_{\alpha }^{T}(x)=\mu_{A}(x)+\alpha$ and $(\nu_{A})_{\alpha }^{T}(x)=\nu_{A}(x)-\alpha$ for all $x\in X$.
\end{definition}

\begin{definition}
Let $A=(\mu_{A},\nu_{A}) $ be an intuitionistic fuzzy subset of a set $X,$ $\beta\in \lbrack 0,1].$ An object having the form $A^{M}_{\beta}=((\mu_{A})^{M}_{\beta},(\nu_{A})^{M}_{\beta})$ is called an {\it intuitionistic fuzzy multiplication} of $A$ if $(\mu_{A})_{\beta}^{M}(x)=\beta\cdot\mu_{A}(x)$ and $(\nu_{A})_{\beta}^{M}(x)=\beta\cdot\nu_{A}(x)$ for all $x\in X$.
\end{definition}

\begin{definition}
Let $A=(\mu_{A},\nu_{A}) $ be an intuitionistic fuzzy subset of a set $X,$ $\alpha \in\lbrack 0,\inf\{\beta\cdot\nu_{A}(x):x\in X\}],$ where $\beta\in \lbrack 0,1].$ An object having the form $A^{C}_{\beta\alpha}=((\mu_{A})^{C}_{\beta\alpha},(\nu_{A})^{C}_{\beta\alpha})$ is called an {\it intuitionistic fuzzy magnified translation} of $A$ if $(\mu_{A})_{\beta\alpha }^{C}(x)=\beta\cdot\mu_{A}(x)+\alpha$ and $(\nu_{A})_{\beta\alpha }^{C}(x)=\beta\cdot\nu_{A}(x)-\alpha$ for
all $x\in X$.
\end{definition}

\begin{example} Let $X=\{1,\omega,\omega^{2}\}.$ Let $A=(\mu_{A},\nu_{A}) $ be an intuitionistic fuzzy subset of $X,$ defined
as follows%
\begin{align*}
\mu_{A} (x)=\left\{
\begin{array}{ll}
0.3 & \text{if} \ x=1 \\
0.1 & \text{if} \ x=\omega \\
0.5 & \text{if} \ x=\omega^{2}
\end{array}
\right. & \text{ and }
\nu_{A} (x)=\left\{
\begin{array}{ll}
0.4 & \text{if} \ x=1 \\
0.25 & \text{if} \ x=\omega \\
0.3 & \text{if} \ x=\omega^{2}
\end{array}
\right..
\end{align*}

Let $\alpha=0.04$ and $\beta=0.2.$ Then the intuitionistic fuzzy magnified translation of $A$ is given by
\begin{align*}
(\mu_{A})_{\beta\alpha }^{C}(x)=\left\{
\begin{array}{ll}
0.1 & \text{if} \ x=1 \\
0.06 & \text{if} \ x=\omega \\
0.14 & \text{if} \ x=\omega^{2}
\end{array}
\right. & \text{ and }
(\nu_{A})_{\beta\alpha }^{C}(x)=\left\{
\begin{array}{ll}
0.04 & \text{if} \ x=1 \\
0.01 & \text{if} \ x=\omega \\
0.02 & \text{if} \ x=\omega^{2}
\end{array}
\right..
\end{align*}
\end{example}

\indent In what follows unless otherwise mentioned $A=(\mu_{A},\nu_{A})$ denotes a non-empty intuitionistic fuzzy subset of $S$ and $A_{\beta\alpha}^{C}$ denotes the intuitionistic fuzzy magnified translation of $A$ where $\beta\in(0,1],\alpha\in\lbrack0,\inf\{\beta\cdot\nu_{A}(x):x\in$ Supp$(\nu_{A})\}].$ It can be noted here that
$A_{\beta\alpha}^{C}$ is also non-empty.

\begin{theorem}
The intuitionistic fuzzy magnified translation $A_{\beta\alpha}^{C}=((\mu_{A})_{\beta\alpha}^{C},(\nu_{A})_{\beta\alpha}^{C})$ of $A=(\mu_{A},\nu_{A})$ is an intuitionistic fuzzy subsemigroup of $S$ if and only if $A=(\mu_{A},\nu_{A})$ is an intuitionistic fuzzy subsemigroup of $S.$
\end{theorem}

\begin{proof}
Let $A=(\mu_{A},\nu_{A})$ be an intuitionitic fuzzy subsemigroup of a semigroup $S.$ Then $A$ is a non-empty intuitionistic fuzzy subset of $S(\beta>0).$ Hence $A_{\beta\alpha}^{C}$ is also non-empty. Now for $x,y\in S$,%

\begin{align*}
(\mu_{A})_{\beta\alpha}^{C}(xy)  &  =\beta\cdot\mu_{A}(xy)+\alpha\\
&  \geq\beta\cdot\min\{\mu_{A}(x),\mu_{A}(y)\}+\alpha(\text{since }A\text{ is an intuitionistic
fuzzy}\\
& \text{ subsemigroup of }S)=\min\{\beta\cdot\mu_{A}(x)+\alpha,\beta\cdot\mu_{A}(y)+\alpha\}\\
&  =\min\{(\mu_{A})_{\beta\alpha}^{C}(x),(\mu_{A})_{\beta\alpha}^{C}(y)\}
\end{align*} and%
\begin{align*}
(\nu_{A})_{\beta\alpha}^{C}(xy)  &  =\beta\cdot\nu_{A}(xy)-\alpha\\
&  \leq\beta\cdot\max\{\nu_{A}(x),\nu_{A}(y)\}-\alpha(\text{since }A\text{ is an intuitionistic
fuzzy}\\
& \text{ subsemigroup of }S)=\max\{\beta\cdot\nu_{A}(x)-\alpha,\beta\cdot\nu_{A}(y)-\alpha\}\\
&  =\max\{(\nu_{A})_{\beta\alpha}^{C}(x),(\nu_{A})_{\beta\alpha}^{C}(y)\}.
\end{align*}

Hence $A_{\beta\alpha}^{C}=((\mu_{A})_{\beta\alpha}^{C},(\nu_{A})_{\beta\alpha}^{C})$ is an intuitionistic fuzzy subsemigroup of $S.$

Conversely, let $A_{\beta\alpha}^{C}=((\mu_{A})_{\beta\alpha}^{C},(\nu_{A})_{\beta\alpha}^{C})$ be an intuitionistic fuzzy subsemigroup of $S.$ Then $A_{\beta\alpha}^{C}$ and hence $A$ is a non-empty fuzzy subset of $S.$ Now for all $x,y\in S,$%
\begin{align*}
(\mu_{A})_{\beta\alpha}^{C}(xy)  &  \geq\min\{(\mu_{A})_{\beta\alpha}^{C}(x),(\mu
_{A})_{\beta\alpha}^{C}(y)\}\\
i.e.,\ \beta.\mu_{A}(xy)+\alpha &  \geq\min\{\beta.\mu_{A}(x)+\alpha,\beta
.\mu_{A}(y)+\alpha\}\\
i.e.,\ \beta.\mu_{A}(xy)+\alpha &  \geq\beta.\min\{\mu_{A}(x),\mu_{A}(y)\}+\alpha\\
i.e.,\ \mu_{A}(xy)  &  \geq\min\{\mu_{A}(x),\mu_{A}(y)\}
\end{align*} and%
\begin{align*}
(\nu_{A})_{\beta\alpha}^{C}(xy)  &  \leq\max\{(\nu_{A})_{\beta\alpha}^{C}(x),(\nu
_{A})_{\beta\alpha}^{C}(y)\}\\
i.e.,\ \beta.\mu_{A}(xy)-\alpha &  \leq\max\{\beta.\nu_{A}(x)-\alpha,\beta
.\nu_{A}(y)-\alpha\}\\
i.e.,\ \beta.\nu_{A}(xy)-\alpha &  \leq\beta.\max\{\nu_{A}(x),\nu_{A}(y)\}-\alpha\\
i.e.,\ \nu_{A}(xy)  &  \leq\max\{\nu_{A}(x),\nu_{A}(y)\}
\end{align*}

Hence $A=(\mu_{A},\nu_{A})$ is an intuitionistic fuzzy subsemigroup of $S.$
\end{proof}

\begin{theorem}
The intuitionistic fuzzy magnified translation $A_{\beta\alpha}^{C}=((\mu_{A})_{\beta\alpha}^{C},(\nu_{A})_{\beta\alpha}^{C})$ of $A=(\mu_{A},\nu_{A})$ is an intuitionistic fuzzy bi-ideal of $S$ if and only if $A=(\mu_{A},\nu_{A})$ is an intuitionistic fuzzy bi-ideal of $S.$
\end{theorem}

\begin{proof}
Let $A=(\mu_{A},\nu_{A})$ be an intuitionistic fuzzy bi-ideal of a semigroup $S.$ Then $A$ is an intuitionistic fuzzy
subsemigroup of $S$ and hence, by Theorem $3.4,$ $A_{\beta\alpha}^{C}$ is an intuitionistic fuzzy subsemigroup of $S.$ Now for $x,\omega,y\in S,$%
\begin{align*}
(\mu_{A})_{\beta\alpha}^{C}(x\omega y)=\beta\cdot\mu_{A}(x\omega y)+\alpha &  \geq\beta\cdot\min\{\mu_{A}(x),\mu_{A}(y)\}+\alpha\\&(\text{since }A\text{ is an intuitionistic
fuzzy}\\
& \text{ bi-ideal of }S)\\
&  =\min\{\beta\cdot\mu_{A}(x)+\alpha,\beta\cdot\mu_{A}(y)+\alpha\}\\
&  =\min\{(\mu_{A})_{\beta\alpha}^{C}(x),(\mu_{A})_{\beta\alpha}^{C}(y)\}.
\end{align*} and%
\begin{align*}
(\nu_{A})_{\beta\alpha}^{C}(x\omega y)=\beta\cdot\nu_{A}(x\omega y)-\alpha &  \leq\beta\cdot\max\{\nu_{A}(x),\nu_{A}(y)\}-\alpha\\&(\text{since }A\text{ is an intuitionistic
fuzzy}\\
& \text{ bi-ideal of }S)\\
& =\max\{\beta\cdot\nu_{A}(x)-\alpha,\beta\cdot\nu_{A}(y)-\alpha\}\\
&  =\max\{(\nu_{A})_{\beta\alpha}^{C}(x),(\nu_{A})_{\beta\alpha}^{C}(y)\}.
\end{align*}
Hence $A_{\beta\alpha}^{C}=((\mu_{A})_{\beta\alpha}^{C},(\nu_{A})_{\beta\alpha}^{C})$ is an intuitionistic fuzzy bi-ideal of $S.$

Conversely, let $A_{\beta\alpha}^{C}=((\mu_{A})_{\beta\alpha}^{C},(\nu_{A})_{\beta\alpha}^{C})$ be an intuitionistic fuzzy bi-ideal of $S.$ Then $A_{\beta\alpha}^{C}=((\mu_{A})_{\beta\alpha}^{C},(\nu_{A})_{\beta\alpha}^{C})$ is an intuitionistic fuzzy subsemigroup of $S,$ and hence by Theorem $3.4,$ $A=(\mu_{A},\nu_{A})$ is an intuitionistic fuzzy subsemigroup of $S.$
Now for all $x,\omega,y\in S,$%
\begin{align*}
(\mu_{A})_{\beta\alpha}^{C}(x\omega y)  &  \geq\min\{(\mu_{A})_{\beta\alpha}^{C}(x),(\mu
_{A})_{\beta\alpha}^{C}(y)\}\\
i.e.,\ \beta.\mu_{A}(x\omega y)+\alpha &  \geq\min\{\beta.\mu_{A}(x)+\alpha,\beta
.\mu_{A}(y)+\alpha\}\\
i.e.,\ \beta.\mu_{A}(x\omega y)+\alpha &  \geq\beta.\min\{\mu_{A}(x),\mu_{A}(y)\}+\alpha\\
i.e.,\ \mu_{A}(x\omega y)  &  \geq\min\{\mu_{A}(x),\mu_{A}(y)\}
\end{align*} and%
\begin{align*}
(\nu_{A})_{\beta\alpha}^{C}(x\omega y)  &  \leq\max\{(\nu_{A})_{\beta\alpha}^{C}(x),(\nu
_{A})_{\beta\alpha}^{C}(y)\}\\
i.e.,\ \beta.\nu_{A}(x\omega y)-\alpha &  \leq\max\{\beta.\nu_{A}(x)-\alpha,\beta
.\nu_{A}(y)-\alpha\}\\
i.e.,\ \beta.\nu_{A}(x\omega y)-\alpha &  \leq\beta.\max\{\nu_{A}(x),\nu_{A}(y)\}-\alpha\\
i.e.,\ \nu_{A}(x\omega y)  &  \leq\max\{\nu_{A}(x),\nu_{A}(y)\}
\end{align*}

Hence $A=(\mu_{A},\nu_{A})$ is an intuitionistic fuzzy bi-ideal of $S.$
\end{proof}

\begin{theorem}
Let $A=(\mu_{A},\nu_{A})$ be a non-empty intuitionistic fuzzy subset of a semigroup $S.$ Then
$A=(\mu_{A},\nu_{A})$ is an intuitionistic fuzzy bi-ideal of $S$ if and only if the intuitionistic fuzzy magnified
translation $A_{\beta\alpha}^{C}=((\mu_{A})_{\beta\alpha}^{C},(\nu_{A})_{\beta\alpha}^{C})$ of $A$ is a constant function, provided $S$ is a group with identity $e.$
\end{theorem}

\begin{proof}
Let $A=(\mu_{A},\nu_{A})$ be an intuitionistic fuzzy bi-ideal of $S.$ Then%
\begin{align*}
(\mu_{A})_{\beta\alpha}^{C}(x)  &  =\beta\cdot\mu_{A}(x)+\alpha=\beta\cdot
\mu_{A}(exe)+\alpha(\text{since }e\text{ is the identity of }S)\\
&  \geq\beta\cdot\min\{\mu_{A}(e),\mu_{A}(e)\}+\alpha(\text{since }A\text{ is an intuitionistic
fuzzy}\\
& \text{ bi-ideal of }S)=\beta\cdot\mu_{A}(e)+\alpha=(\mu_{A})_{\beta\alpha}^{C}(e)
\end{align*}
Again%
\begin{align*}
(\mu_{A})_{\beta\alpha}^{C}(e)  &  =(\mu_{A})_{\beta\alpha}^{C}(ee)=\beta\cdot
\mu_{A}(ee)+\alpha(\text{since }e\text{ is the identity of }S)\\
&  =\beta\cdot\mu_{A}((xx^{-1})(x^{-1}x))+\alpha\\
&  =\beta\cdot\mu_{A}(x(x^{-1}x^{-1})x)+\alpha\\
&  \geq\beta\cdot\min\{\mu_{A}(x),\mu_{A}(x)\}+\alpha(\text{since }A\text{ is an intuitionistic
fuzzy}\\
& \text{ bi-ideal of }S)=\beta\cdot\mu_{A}(x)+\alpha=(\mu_{A})_{\beta\alpha}^{C}(x).
\end{align*}
Thus $(\mu_{A})_{\beta\alpha}^{C}(x)=(\mu_{A})_{\beta\alpha}^{C}(e)$ $\forall x\in S.$ Similarly we can show that $(\nu_{A})_{\beta\alpha}^{C}(x)=(\nu_{A})_{\beta\alpha}^{C}(e)$ $\forall x\in S.$ Hence $A_{\beta\alpha}^{C}$ is a constant function.

Conversely, suppose $A_{\beta\alpha}^{C}$ is a constant function. Then $A_{\beta\alpha}^{C}$ is an intuitionistic fuzzy bi-ideal of $S\cite{K1}.$ Hence by Theorem $3.5,$ $A$ is an intuitionistic fuzzy bi-ideal of $S.$
\end{proof}

\begin{note}
For the converse to be true, $S$ need not be a group.
\end{note}

\begin{theorem}
The intuitionistic fuzzy magnified translation $A_{\beta\alpha}^{C}=((\mu_{A})_{\beta\alpha}^{C},(\nu_{A})_{\beta\alpha}^{C})$ of $A=(\mu_{A},\nu_{A})$ is an intuitionistic fuzzy $(1,2)$-ideal of $S$ if and only if $A=(\mu_{A},\nu_{A})$ is an intuitionistic fuzzy $(1,2)$-ideal of $S.$
\end{theorem}

\begin{proof}
Let $A=(\mu_{A},\nu_{A})$ be an intuitionistic fuzzy $(1,2)$-ideal of $S.$ Then $A$ is an intuitionistic fuzzy subsemigroup of $S.$ Hence, by Theorem $3.4,$ $A_{\beta\alpha}^{C}$ is an intuitionistic fuzzy subsemigroup of $S.$ Now for $x,y,z,\omega\in S,$%
\begin{align*}
(\mu_{A})_{\beta\alpha}^{C}(x\omega(yz))  &  =\beta\cdot\mu_{A}(x\omega(yz))+\alpha\\
&  \geq\beta\cdot\min\{\mu_{A}(x),\mu_{A}(y),\mu_{A}(z)\}+\alpha\\
&  (\text{since }A\text{ is an intuitionistic fuzzy }(1,2)\text{-ideal of }S)\\
&  =\min\{\beta\cdot\mu_{A}(x)+\alpha,\beta\cdot\mu_{A}(y)+\alpha,\beta\cdot
\mu_{A}(z)+\alpha\}\\
&  =\min\{(\mu_{A})_{\beta\alpha}^{C}(x),(\mu_{A})_{\beta\alpha}^{C}(y),(\mu
_{A})_{\beta\alpha}^{C}(z)\}.
\end{align*} and%
\begin{align*}
(\nu_{A})_{\beta\alpha}^{C}(x\omega(yz))  &  =\beta\cdot\nu_{A}(x\omega(yz))-\alpha\\
&  \leq\beta\cdot\max\{\nu_{A}(x),\nu_{A}(y),\nu_{A}(z)\}-\alpha\\
&  (\text{since }A\text{ is an intuitionistic fuzzy }(1,2)\text{-ideal of }S)\\
&  =\max\{\beta\cdot\nu_{A}(x)-\alpha,\beta\cdot\nu_{A}(y)-\alpha,\beta\cdot
\nu_{A}(z)-\alpha\}\\
&  =\max\{(\nu_{A})_{\beta\alpha}^{C}(x),(\nu_{A})_{\beta\alpha}^{C}(y),(\nu
_{A})_{\beta\alpha}^{C}(z)\}.
\end{align*}
Hence $A_{\beta\alpha}^{C}=((\mu_{A})_{\beta\alpha}^{C},(\nu_{A})_{\beta\alpha}^{C})$ is an intuitionistic fuzzy $(1,2)$-ideal of $S.$

Conversely, let $A_{\beta\alpha}^{C}=((\mu_{A})_{\beta\alpha}^{C},(\nu_{A})_{\beta\alpha}^{C})$ be an intuitionistic fuzzy $(1,2)$-ideal of $S.$ Then $A_{\beta\alpha}^{C}=((\mu_{A})_{\beta\alpha}^{C},(\nu_{A})_{\beta\alpha}^{C})$ is an intuitionistic fuzzy subsemigroup of $S,$ and hence by Theorem $3.4,$ $A=(\mu_{A},\nu_{A})$ is an intuitionistic fuzzy subsemigroup of $S.$
Now for all $x,\omega,y,z\in S,$%
\begin{align*}
(\mu_{A})_{\beta\alpha}^{C}(x\omega (yz))  &  \geq\min\{(\mu_{A})_{\beta\alpha}^{C}(x),(\mu
_{A})_{\beta\alpha}^{C}(y),(\mu_{A})_{\beta\alpha}^{C}(z)\}\\
i.e.,\ \beta.\mu_{A}(x\omega (yz))+\alpha &  \geq\min\{\beta.\mu_{A}(x)+\alpha,\beta
.\mu_{A}(y)+\alpha,\beta.\mu_{A}(z)+\alpha\}\\
i.e.,\ \beta.\mu_{A}(x\omega (yz))+\alpha &  \geq\beta.\min\{\mu_{A}(x),\mu_{A}(y),\mu_{A}(z)\}+\alpha\\
i.e.,\ \mu_{A}(x\omega y)  &  \geq\min\{\mu_{A}(x),\mu_{A}(y),\mu_{A}(z)\}
\end{align*} and%
\begin{align*}
(\nu_{A})_{\beta\alpha}^{C}(x\omega (yz))  &  \leq\max\{(\nu_{A})_{\beta\alpha}^{C}(x),(\nu
_{A})_{\beta\alpha}^{C}(y),(\nu_{A})_{\beta\alpha}^{C}(z)\}\\
i.e.,\ \beta.\nu_{A}(x\omega (yz))-\alpha &  \leq\max\{\beta.\nu_{A}(x)-\alpha,\beta
.\nu_{A}(y)-\alpha,\beta.\nu_{A}(z)-\alpha\}\\
i.e.,\ \beta.\nu_{A}(x\omega (yz))-\alpha &  \leq\beta.\max\{\nu_{A}(x),\nu_{A}(y),\nu_{A}(z)\}-\alpha\\
i.e.,\ \nu_{A}(x\omega (yz))  &  \leq\max\{\nu_{A}(x),\nu_{A}(y),\nu_{A}(z)\}
\end{align*}

Hence $A=(\mu_{A},\nu_{A})$ is an intuitionistic fuzzy $(1,2)$-ideal of $S.$
\end{proof}

\begin{theorem}
The intuitionistic fuzzy magnified translation $A_{\beta\alpha}^{C}=((\mu_{A})_{\beta\alpha}^{C},(\nu_{A})_{\beta\alpha}^{C})$ of $A=(\mu_{A},\nu_{A})$ is an intuitionistic fuzzy left ideal$($intuitionistic fuzzy right ideal, intuitionistic fuzzy ideal$)$ of $S$ if and only if $A=(\mu_{A},\nu_{A})$ is an intuitionstic fuzzy left ideal$($resp. intuitionistic fuzzy right ideal, intuitionstic fuzzy ideal$)$ of $S.$
\end{theorem}

\begin{proof}
Let $A=(\mu_{A},\nu_{A})$ be an intuitionistic fuzzy left ideal of $S.$ Then $A=(\mu_{A},\nu_{A})$ is a non-empty intuitionistic fuzzy subset of $S.$ Hence, as $\beta>0,$ $A_{\beta\alpha}^{C}$ is a non-empty intuitionistic fuzzy
subset of $S.$ Let $x,y\in S.$ Then%
\[
(\mu_{A})_{\beta\alpha}^{C}(xy)=\beta\cdot\mu_{A}(xy)+\alpha\geq\beta\cdot\mu_{A}
(y)+\alpha=(\mu_{A})_{\beta\alpha}^{C}(y)
\] and%
\[
(\nu_{A})_{\beta\alpha}^{C}(xy)=\beta\cdot\nu_{A}(xy)-\alpha\leq\beta\cdot\nu_{A}
(y)-\alpha=(\nu_{A})_{\beta\alpha}^{C}(y).
\]
Hence $A_{\beta\alpha}^{C}=((\mu_{A})_{\beta\alpha}^{C},(\nu_{A})_{\beta\alpha}^{C})$ is an intuitionistic fuzzy left ideal $S.$

Conversely, let $A_{\beta\alpha}^{C}=((\mu_{A})_{\beta\alpha}^{C},(\nu_{A})_{\beta\alpha}^{C})$ be an intuitionistic fuzzy left ideal of $S.$ Then $A_{\beta\alpha}^{C}$ and hence $A$ is a non-empty intuitionistic fuzzy subset of $S.$
Now for all $x,y\in S,$%
\[
(\mu_{A})_{\beta\alpha}^{C}(xy)\geq(\mu_{A})_{\beta\alpha}^{C}(y),i.e.,\ \beta
.\mu_{A}(xy)+\alpha\geq\beta.\mu_{A}(y)+\alpha,i.e.,\ \mu_{A}(xy)\geq\mu_{A}(y)
\]and%
\[
(\nu_{A})_{\beta\alpha}^{C}(xy)\leq(\nu_{A})_{\beta\alpha}^{C}(y),i.e.,\ \beta
.\nu_{A}(xy)-\alpha\leq\beta.\nu_{A}(y)-\alpha,i.e.,\ \nu_{A}(xy)\leq\nu_{A}(y).
\]
Hence $A=(\mu_{A},\nu_{A})$ is an intuitionistic fuzzy left ideal of $S.$ Similar is the proof for intuitionistic fuzzy right ideal or intuitionistic fuzzy ideal.
\end{proof}

\begin{theorem}
The intuitionistic fuzzy magnified translation $A_{\beta\alpha}^{C}=((\mu_{A})_{\beta\alpha}^{C},(\nu_{A})_{\beta\alpha}^{C})$ of $A=(\mu_{A},\nu_{A})$ is an intuitionistic fuzzy semiprime ideal of $S$ if and only if $A=(\mu_{A},\nu_{A})$ is an intuitionistic fuzzy semiprime ideal of $S.$
\end{theorem}

\begin{proof}
Let $A=(\mu_{A},\nu_{A})$ be an intuitionistic fuzzy semiprime ideal of $S.$ Then $A=(\mu_{A},\nu_{A})$ is a non-empty intuitionistic fuzzy subset of $S.$ Hence $A_{\beta\alpha}^{C}$ is a non-empty intuitionistic fuzzy subset of
$S(\beta>0).$ Let $x\in S.$ Then%
\[
(\mu_{A})_{\beta\alpha}^{C}(x)=\beta\cdot\mu_{A}(x)+\alpha\geq\beta\cdot\mu_{A}
(x^{2})+\alpha=(\mu_{A})_{\beta\alpha}^{C}(x^{2})
\] and%
\[
(\nu_{A})_{\beta\alpha}^{C}(x)=\beta\cdot\nu_{A}(x)-\alpha\leq\beta\cdot\nu_{A}
(x^{2})-\alpha=(\nu_{A})_{\beta\alpha}^{C}(x^{2}).
\]
Hence $A_{\beta\alpha}^{C}=((\mu_{A})_{\beta\alpha}^{C},(\nu_{A})_{\beta\alpha}^{C})$ is an intuitionistic fuzzy semiprime ideal of $S.$

Conversely, let $A_{\beta\alpha}^{C}=((\mu_{A})_{\beta\alpha}^{C},(\nu_{A})_{\beta\alpha}^{C})$ be an intuitionistic fuzzy semiprime ideal of $S.$ Then $A_{\beta\alpha}^{C}$ and hence $A$ is a non-empty intuitionistic fuzzy subset of
$S.$ Then for all $x\in S,$%
\[
(\mu_{A})_{\beta\alpha}^{C}(x)\geq(\mu_{A})_{\beta\alpha}^{C}(x^{2}),i.e.,\ \beta
.\mu_{A}(x)+\alpha\geq\beta.\mu_{A}(x^{2})+\alpha,i.e.,\ \mu_{A}(x)\geq\mu_{A}(x^{2}).
\] and%
\[
(\nu_{A})_{\beta\alpha}^{C}(x)\leq(\nu_{A})_{\beta\alpha}^{C}(x^{2}),i.e.,\ \beta
.\nu_{A}(x)-\alpha\leq\beta.\nu_{A}(x^{2})-\alpha,i.e.,\ \nu_{A}(x)\leq\nu_{A}(x^{2}).
\]
Hence $A=(\mu_{A},\nu_{A})$ is an intuitionistic fuzzy semiprime ideal of $S.$
\end{proof}

\begin{theorem}
Let $A=(\mu_{A},\nu_{A})$ and $B=(\mu_{B},\nu_{B})$ be two intuitionistic fuzzy semiprime ideals of $S.$ Then $A\cap B$ is an intuitionistic fuzzy semiprime ideal of $S,$ provided it is non-empty$.$
\end{theorem}

\begin{proof}
Let $A=(\mu_{A},\nu_{A})$ and $B=(\mu_{B},\nu_{B})$ be two intuitionistic fuzzy semiprime ideals of $S$ and $x\in S.$ Then%
\begin{align*}
(\mu_{A}\cap\mu_{B})(x)  &  = \min\{\mu_{A}(x),\mu_{B}(x)\}\geq\min\{\mu_{A}(x^{2}),\mu_{B}(x^{2})\}\\
&  (\text{since }A \text{ and } B\text{ are intuitionistic fuzzy semiprime ideals of }S)\\
&  =(\mu_{A}\cap\mu_{B})(x^{2})
\end{align*} and%
\begin{align*}
(\nu_{A}\cup\nu_{B})(x)  &  = \max\{\nu_{A}(x),\nu_{B}(x)\}\leq\max\{\nu_{A}(x^{2}),\nu_{B}(x^{2})\}\\
&  (\text{since }A \text{ and } B\text{ are intuitionistic fuzzy semiprime ideals of }S)\\
&  =(\nu_{A}\cup\nu_{B})(x^{2}).
\end{align*}

Hence $A\cap B$ is an intuitionistic fuzzy semiprime ideal of $S.$
\end{proof}

Since the intersections of intuitionistic fuzzy semiprime ideals is an intuitionistic fuzzy semiprime ideal, so we have the following corollary.

\begin{corollary}
Let $A=(\mu_{A},\nu_{A})$ and $B=(\mu_{B},\nu_{B})$ be two intuitionistic fuzzy semiprime ideals of $S.$ Then $A_{\beta\alpha}^{C}\cap B_{\beta\alpha}^{C}$ is an intuitionistic fuzzy semiprime ideal of $S,$ provided it is non-empty$.$
\end{corollary}

\begin{notation}
Let $A=(\mu_{A},\nu_{A})$ be an intuitionistic fuzzy subset of a semigroup $S$ and $A_{\beta\alpha}^{C}=((\mu_{A})_{\beta\alpha}^{C},(\nu_{A})_{\beta\alpha}^{C})$ be an intuitionistic fuzzy magnified translation of $A.$ If for any elements $x,y\in S,$ $(\mu_{A})_{\beta\alpha}^{C}(x)=(\mu_{A})_{\beta\alpha}^{C}(y)$ and $(\nu_{A})_{\beta\alpha}^{C}(x)=(\nu_{A})_{\beta\alpha}^{C}(y),$ then we write $A_{\beta\alpha}^{C}(x)=A_{\beta\alpha}^{C}(y)$ $\forall x,y\in S.$
\end{notation}

\begin{theorem}
If $A=(\mu_{A},\nu_{A})$ is an intuitionistic fuzzy semiprime ideal of a semigroup $S,$ then $A_{\beta\alpha}^{C}(x)=A_{\beta\alpha}^{C}(x^{2})$ $\forall x\in S.$
\end{theorem}

\begin{proof}
Let $A=(\mu_{A},\nu_{A})$ be an intuitionistic fuzzy semiprime ideal of $S.$ Let $x\in S.$ Then%
\[
(\mu_{A})_{\beta\alpha}^{C}(x)=\beta\cdot\mu_{A}(x)+\alpha\geq\beta\cdot\mu_{A}
(x^{2})+\alpha=(\mu_{A})_{\beta\alpha}^{C}(x^{2}).
\]
Again
\begin{align*}
(\mu_{A})_{\beta\alpha}^{C}(x^{2})  &  =\beta\cdot\mu_{A}(x^{2})+\alpha\\
&  \geq\beta\cdot\mu_{A}(x)+\alpha(\text{since }A\text{ is an intuitionistic fuzzy ideal})\\
&  =(\mu_{A})_{\beta\alpha}^{C}(x).
\end{align*}
 Consequently, $(\mu_{A})_{\beta\alpha}^{C}(x)=(\mu_{A})_{\beta\alpha}^{C}(x^{2})$ $\forall x\in S.$ Similarly we can show that $(\nu_{A})_{\beta\alpha}^{C}(x)=(\nu_{A})_{\beta\alpha}^{C}(x^{2})$ $\forall x\in S.$ Hence $A_{\beta\alpha}^{C}(x)=A_{\beta\alpha}^{C}(x^{2})$ $\forall x\in S.$
\end{proof}

By routine verification we can have the following theorem.

\begin{theorem}
Let $\chi$ be the characteristic function of a non-empty subset $A$ of $S.$ Then $A$ is a left$($resp. right$)$ ideal of $S$ if and only if $(\chi,\overline{\chi})$ is an intuitionistic fuzzy left$($resp. right$)$ ideal of $S.$
\end{theorem}

\begin{definition}
A semigroup $S$ is called {\it intra-regular}\cite{K2} if for each element $a$ of $S,$ there exist elements $x,y\in S$ such that $a=xa^{2}y.$
\end{definition}

\begin{theorem}
For a semigroup $S$ the following conditions are equivalent:

$(1)$ $S$ is an intra-regular semigroup,

$(2)$ for every intuitionistic fuzzy ideal $A=(\mu_{A},\nu_{A})$ of $S$ the intuitionistic fuzzy magnified translation $A_{\beta\alpha}^{C}=((\mu_{A})_{\beta\alpha}^{C},(\nu_{A})_{\beta\alpha}^{C})$ of $A$ is an intuitionistic fuzzy semiprime ideal of $S$.
\end{theorem}

\begin{proof}
$(1)\Rightarrow(2):$ Let $A=(\mu_{A},\nu_{A})$ be an intuitionistic fuzzy ideal of $S$ and $m\in$ $S.$ Then there exist $x,y\in S$ such that $m=xm^{2}y($since $S$ is intra-regular$).$ Then%
\begin{align*}
(\mu_{A})_{\beta\alpha}^{C}(m)  &  =(\mu_{A})_{\beta\alpha}^{C}(xm^{2}y)=\beta\cdot
\mu_{A}(xm^{2}y)+\alpha\\
&  \geq\beta\cdot\mu_{A}(m^{2}y)+\alpha\geq\beta\cdot\mu_{A}(m^{2})+\alpha=(\mu
_{A})_{\beta\alpha}^{C}(m^{2})
\end{align*} and
\begin{align*}
(\nu_{A})_{\beta\alpha}^{C}(m)  &  =(\nu_{A})_{\beta\alpha}^{C}(xm^{2}y)=\beta\cdot
\nu_{A}(xm^{2}y)-\alpha\\
&  \leq\beta\cdot\nu_{A}(m^{2}y)-\alpha\leq\beta\cdot\nu_{A}(m^{2})-\alpha=(\nu
_{A})_{\beta\alpha}^{C}(m^{2})
\end{align*}
Hence $A_{\beta\alpha}^{C}=((\mu_{A})_{\beta\alpha}^{C},(\nu_{A})_{\beta\alpha}^{C})$ is an intuitionistic fuzzy semiprime ideal of $S.$

$(2)\Rightarrow(1):$ Let $m$ be any element of $S.$ Then it follows that the intuitionistic fuzzy subset $A=(\chi_{<m^{2}>},\overline{\chi}_{<m^{2}>})$ of the principal ideal $<m^{2}>$ of $S($generated by $m^{2})$ is an intuitionistic fuzzy ideal of $S($ where $\chi_{<m^{2}>}$ is the characteristic function of the principal ideal $<m^{2}>$ of $S).$ From $(2),$  $A_{\beta\alpha}^{C}$ is an intuitionistic fuzzy semiprime ideal of $S.$ Hence $A=(\mu_{A},\nu_{A})$ is an intuitionistic fuzzy semiprime ideal of $S(cf.$ Theorem $3.10).$ So $A_{\beta\alpha}^{C}(m)=A_{\beta\alpha}^{C}(m^{2})(cf.$ Theorem $3.13).$ Hence
$\chi_{<m^{2}>\text{ }}(m)=\chi_{<m^{2}>\text{ }}(m^{2})$ and $\overline{\chi}_{<m^{2}>\text{ }}(m)=\overline{\chi}_{<m^{2}>\text{ }}(m^{2}).$ Since $m^{2}%
\in<m^{2}>,$ we have $\chi_{<m^{2}>\text{ }}(m)=\chi_{<m^{2}>\text{ }}%
(m^{2})=1$ and $\overline{\chi}_{<m^{2}>\text{ }}(m)=\overline{\chi}_{<m^{2}>\text{ }}%
(m^{2})=0.$ So $m\in<m^{2}>=\{m^{2}\}\cup Sm^{2}\cup m^{2}S\cup Sm^{2}S.$ This
shows that $S$ is intra-regular. This completes the proof.
\end{proof}

\begin{definition}
A semigroup $S$ is said to be {\it left $($right$)$ regular}\cite{K1} if, for each element $a$ of $S$, there exists an element $x$ in $S$ such that $a=xa^{2}($resp. $a=a^{2}x).$
\end{definition}

\begin{theorem}
For a semigroup $S$ the following conditions are equivalent:

$(1)$ $S$ is a left regular semigroup,

$(2)$ for every intuitionistic fuzzy left ideal $A=(\mu_{A},\nu_{A})$ of $S$ the intuitionistic fuzzy magnified translation $A_{\beta\alpha}^{C}=((\mu_{A})_{\beta\alpha}^{C},(\nu_{A})_{\beta\alpha}^{C})$ of $A=(\mu_{A},\nu_{A})$ is an intuitionistic fuzzy semiprime ideal of $S$.
\end{theorem}

\begin{proof}
$(1)\Rightarrow(2):$ Let $A=(\mu_{A},\nu_{A})$  be an intuitionistic fuzzy left ideal of $S$. Then by Theorem $3.11,$ $A_{\beta\alpha}^{C}$ is an intuitionistic fuzzy left ideal of $S.$ Let $m\in S$. Then there exists an element $x\in
S$ such that $m=xm^{2}($since $S$ is left regular$)$. Then%
\begin{align*}
(\mu_{A})_{\beta\alpha}^{C}(m)  &  =\beta\cdot\mu_{A}(m)+\alpha=\beta\cdot\mu_{A}
(xm^{2})+\alpha\\
&  \geq\beta\cdot\mu_{A}(m^{2})+\alpha(\text{since }A\text{ is an intuitionistic fuzzy left
ideal of }S)\\
&  =(\mu_{A})_{\beta\alpha}^{C}(m^{2})
\end{align*} and%
\begin{align*}
(\nu_{A})_{\beta\alpha}^{C}(m)  &  =\beta\cdot\nu_{A}(m)-\alpha=\beta\cdot\nu_{A}
(xm^{2})-\alpha\\
&  \leq\beta\cdot\nu_{A}(m^{2})-\alpha(\text{since }A\text{ is an intuitionistic fuzzy left
ideal of }S)\\
&  =(\nu_{A})_{\beta\alpha}^{C}(m^{2}).
\end{align*}
Hence $A_{\beta\alpha}^{C}=((\mu_{A})_{\beta\alpha}^{C},(\nu_{A})_{\beta\alpha}^{C})$ is an intuitionistic fuzzy semiprime ideal of $S.$

$(2)\Rightarrow(1):$ Suppose $(2)$ holds. Let $m$ be any element of $S.$ Then it follows that the intuitionistic fuzzy subset $A=(\chi_{<m^{2}|},\overline{\chi}_{<m^{2}|})$ of the left ideal $<m^{2}|$ of $S($generated by $m^{2})$ is an intuitionistic fuzzy left ideal of $S($ where $\chi_{<m^{2}|}$ is the characteristic function of the left ideal $<m^{2}|$ of $S).$ From $(2),$  $A_{\beta\alpha}^{C}$ is an intuitionistic fuzzy semiprime ideal of $S.$ Hence $A=(\mu_{A},\nu_{A})$ is an intuitionistic fuzzy semiprime ideal of $S(cf.$ Theorem $3.10).$ So $A_{\beta\alpha}^{C}(m)=A_{\beta\alpha}^{C}(m^{2})(cf.$ Theorem $3.13).$ Hence $\chi_{<m^{2}|\text{}}(m)=\chi_{<m^{2}|\text{}}(m^{2})$ and $\overline{\chi}_{<m^{2}|\text{}}(m)=\overline{\chi}_{<m^{2}|\text{}}(m^{2}).$ Now since $m^{2}\in<m^{2}|,$ we see that $\chi_{<m^{2}|\text{ }}(m)=\chi_{<m^{2}|\text{ }}(m^{2})=1$ and $\overline{\chi}_{<m^{2}|\text{ }}(m)=\overline{\chi}_{<m^{2}|\text{ }}(m^{2})=0.$ Hence $m\in<m^{2}|=\{m^{2}\}\cup Sm^{2}.$ This shows that $S$ is left regular. This completes the proof.
\end{proof}

In a similar way we can prove the following theorem.

\begin{theorem}
For a semigroup $S$ the following conditions are equivalent:

$(1)$ $S$ is a right regular semigroup,

$(2)$ for every intuitionistic fuzzy right ideal $A=(\mu_{A},\nu_{A})$ of $S$ the intuitionistic fuzzy magnified translation $A_{\beta\alpha}^{C}=((\mu_{A})_{\beta\alpha}^{C},(\nu_{A})_{\beta\alpha}^{C})$ of $A$ is an intuitionistic fuzzy semiprime ideal of $S$.
\end{theorem}

\begin{definition}
A semigroup $S$ is called {\it archimedean}\cite{K2} if for all $a,b\in S,$ there exists a positive integer $n$ such that $a^{n}\in SbS.$
\end{definition}

\begin{theorem}
Let $A=(\mu_{A},\nu_{A})$ be an intuitionistic fuzzy semiprime ideal of an archimedean semigroup $S$.
Then $A_{\beta\alpha}^{C}=((\mu_{A})_{\beta\alpha}^{C},(\nu_{A})_{\beta\alpha}^{C})$ is a constant function.
\end{theorem}

\begin{proof}
Let $m,n\in S.$ Then $S$ being archimedean, there exists a positive integer
$k$ such that $m^{k}=xny$ for some $x,y\in S.$ Now since $A=(\mu_{A},\nu_{A})$ is an intuitionistic fuzzy
semiprime ideal of $S,$ so by Theorem $3.8,$ $A_{\beta\alpha}^{C}$ is an intuitionistic fuzzy semiprime ideal of $S.$ Then%
\begin{align*}
(\mu_{A})_{\beta\alpha}^{C}(m)  &  \geq(\mu_{A})_{\beta\alpha}^{C}(m^{k})=(\mu_{A})_{\beta
\alpha}^{C}(xny)=\beta\cdot\mu_{A}(xny)+\alpha\\
&  \geq\beta\cdot\mu_{A}(n)+\alpha=(\mu_{A})_{\beta\alpha}^{C}(n).
\end{align*}
Using the duality of $m$ and $n$ we deduce that $(\mu_{A})_{\beta\alpha}^{C}%
(n)\geq(\mu_{A})_{\beta\alpha}^{C}(m)$. Thus $(\mu_{A})_{\beta\alpha}^{C}(m)=(\mu
_{A})_{\beta\alpha}^{C}(n)$ $\forall m,n\in S.$ Applying similar argument we can show that $(\nu_{A})_{\beta\alpha}^{C}(m)=(\nu_{A})_{\beta\alpha}^{C}(n)$ $\forall m,n\in S.$ Hence $A_{\beta\alpha}^{C}$ is a constant function.
\end{proof}

\begin{definition}
A semigroup $S$ is called {\it regular}\cite{K2} if for each element $a$ of $S,$ there exists an element $x\in S$ such that $a=axa.$
\end{definition}

\begin{definition}
Let $S$ be a semigroup. Let $A=(\mu_{A},\nu_{A})$ and $B=(\mu_{B},\nu_{B})$ be two intuitionistic fuzzy subsets of $S.$ Then the {\it product} $A\circ B$ of $A$ and $B$ is defined as\\
$$A\circ B=\{<x,(\mu_{A}\circ\mu_{B})(x),(\nu_{A}\circ\nu_{B})(x)>:x\in S\}$$
\end{definition}

$$
\text{ where } (\mu_{A}\circ \mu_{B})(x)=\left\{
\begin{array}
[c]{c}%
\underset{x=uv}{\sup}[\min\{\mu_{A}(u),\mu_{B}(v)\}:u,v\in S]\\
0,\text{ if for any }u,v\in S\text{ },x\neq uv
\end{array}
\right.\\
$$
$$
 \text{ and }
(\nu_{A}\circ\nu_{B})(x)=\left\{
\begin{array}
[c]{c}%
\underset{x=uv}{\inf}[\max\{\nu_{A}(u),\nu_{B}(v)\}:u,v\in S]\\
1,\text{ if for any }u,v\in S\text{ },x\neq uv
\end{array}
\right.
$$

\begin{theorem}
If the semigroup $S$ is both regular and intra-regular then

$(1)$ $A_{\beta\alpha}^{C}$ $\circ B_{\beta\alpha}^{C}\supset A
_{\beta\alpha}^{C}\cap B_{\beta\alpha}^{C}.$

$(2)$ $(A_{\beta\alpha}^{C}\circ B_{\beta\alpha}^{C})\cap
(B_{\beta\alpha}^{C}\circ A_{\beta\alpha}^{C})\supset A_{\beta\alpha
}^{C}\cap B_{\beta\alpha}^{C},$ where $A=(\mu_{A},\nu_{A})$ and $B=(\mu_{B},\nu_{B})$ are intuitionistic fuzzy
bi-ideals of $S$.
\end{theorem}

\begin{proof}
Let $A=(\mu_{A},\nu_{A})$ and $B=(\mu_{B},\nu_{B})$ are intuitionistic fuzzy
bi-ideals of $S$ and $a\in S.$ Then there exist $x,y,z\in S$
such that $a=axa=axaxa($since $S$ is regular$)=ax(ya^{2}z)xa($since $S$ is
intra-regular$)=(axya)(azxa).$ Now
\begin{align*}
\mu_{A}(axya)\geq\min\{\mu_{A}(a),\mu_{A}(a)\}(\text{since }A \text{ is intuitionistic fuzzy bi-ideal of }S)=\mu_{A}(a)
\end{align*}and
\begin{align*}
\nu_{A}(axya)\leq\max\{\nu_{A}(a),\nu_{A}(a)\}(\text{since }A \text{ is intuitionistic fuzzy bi-ideal of }S)=\nu_{A}(a).
\end{align*}
Similarly we can have $\mu_{B}(azxa)\geq\mu_{B}(a)$ and $\nu_{B}(azxa)\leq\nu_{B}(a).$ Then%
\begin{align*}
((\mu_{A})_{\beta\alpha}^{C}\circ(\mu_{B})_{\beta\alpha}^{C})(a)  &  =\underset
{a=pq}{\sup}\min\{(\mu_{A})_{\beta\alpha}^{C}(p),(\mu_{B})_{\beta\alpha}^{C}(q)\}\\
&  =\underset{a=pq}{\sup}\min\{\beta\cdot\mu_{A}(p)+\alpha,\beta\cdot
\mu_{B}(q)+\alpha\}\\
&  =\beta\cdot\underset{a=pq}{\sup}\min\{\mu_{A}(p),\mu_{B}(q)\}+\alpha\\
&  \geq\beta\cdot\min\{\mu_{A}(axya),\mu_{B}(azxa)\}+\alpha(\text{since
}a=(axya)(azxa))\\
&  \geq\beta\cdot\min\{\mu_{A}(a),\mu_{B}(a)\}+\alpha(\text{since }A\text{ and
}B\text{ are intuitionistic}\\
& \text{ fuzzy bi-ideals of }S)=\min\{\beta\cdot\mu_{A}(a)+\alpha,\beta\cdot\mu_{B}(a)+\alpha\}\\
&  =\min\{(\mu_{A})_{\beta\alpha}^{C}(a),\ (\mu_{B})_{\beta\alpha}^{C}(a)\ \}=((\mu
_{A})_{\beta\alpha}^{C}\cap(\mu_{B})_{\beta\alpha}^{C})(a)
\end{align*} again%
\begin{align*}
((\nu_{A})_{\beta\alpha}^{C}\circ(\nu_{B})_{\beta\alpha}^{C})(a)  &  =\underset
{a=pq}{\inf}\max\{(\nu_{A})_{\beta\alpha}^{C}(p),(\nu_{B})_{\beta\alpha}^{C}(q)\}\\
&  =\underset{a=pq}{\inf}\max\{\beta\cdot\nu_{A}(p)-\alpha,\beta\cdot
\nu_{B}(q)-\alpha\}\\
&  =\beta\cdot\underset{a=pq}{\inf}\max\{\nu_{A}(p),\nu_{B}(q)\}-\alpha\\
&  \leq\beta\cdot\max\{\nu_{A}(axya),\nu_{B}(azxa)\}-\alpha(\text{since }a=(axya)(azxa))
\end{align*}
\begin{align*}
&  \leq\beta\cdot\max\{\nu_{A}(a),\nu_{B}(a)\}-\alpha(\text{since }A\text{ and
}B\text{ are intuitionistic fuzzy}\\
& \text{bi-ideals of }S)=\max\{\beta\cdot\nu_{A}(a)-\alpha,\beta\cdot\nu_{B}(a)-\alpha\}\\
&  =\max\{(\nu_{A})_{\beta\alpha}^{C}(a),\ (\nu_{B})_{\beta\alpha}^{C}(a)\ \}=((\nu
_{A})_{\beta\alpha}^{C}\cup(\nu_{B})_{\beta\alpha}^{C})(a)
\end{align*}

Hence $A_{\beta\alpha}^{C}$ $\circ B_{\beta\alpha}^{C}\supset A_{\beta\alpha}^{C}\cap B_{\beta\alpha}^{C}$ and $(1)$ follows. Similarly we can prove that $B_{\beta\alpha}^{C}$ $\circ A_{\beta\alpha}^{C}\supset A_{\beta
\alpha}^{C}\cap B_{\beta\alpha}^{C}$. Combining these two results we can obtain $(2).$ This completes the proof.
\end{proof}

\begin{theorem}
If the semigroup $S$ is both regular, intra-regular and left regular then

$(1)$ $A_{\beta\alpha}^{C}$ $\circ B_{\beta\alpha}^{C}\supset A_{\beta\alpha}^{C}\cap B_{\beta\alpha}^{C}.$

$(2)$ $(A_{\beta\alpha}^{C}\circ B_{\beta\alpha}^{C})\cap(B_{\beta\alpha}^{C}\circ A_{\beta\alpha}^{C})\supset A_{\beta\alpha}^{C}\cap B_{\beta\alpha}^{C}$ where $A=(\mu_{A},\nu_{A})$ and $B=(\mu_{B},\nu_{B})$ are intuitionistic fuzzy $(1,2)$-ideals of $S$.
\end{theorem}

\begin{proof}
Let $A=(\mu_{A},\nu_{A})$ and $B=(\mu_{B},\nu_{B})$ are intuitionistic fuzzy $(1,2)$-ideals of $S$ and $a\in S.$ Then there exist $x,y,z,p\in S$ such that $a=axa=axaxa($since $S$ is regular$)=ax(ya^{2}z)xa($since $S$ is intra-regular$)=(axya)(azxa)=(axypa^{2})(azxpa^{2})(S\,\ $is left regular$).$ Now%
\begin{align*}
\mu_{A}(axypa^{2})  &  =\mu_{A}(axypaa)=\mu_{A}(axyp(aa))\geq\min\{\mu_{A}(a),\mu_{A}(a),\mu_{A}(a)\}(\text{since }A\text{ is }\\
& \text{ intuitionistic fuzzy }(1,2)\text{-ideal of }S)=\mu_{A}(a)
\end{align*} and%
\begin{align*}
\nu_{A}(axypa^{2})  &  =\nu_{A}(axypaa)=\nu_{A}(axyp(aa))\leq\max\{\nu_{A}(a),\nu_{A}(a),\nu_{A}(a)\}(\text{since }A\text{ is }\\
& \text{ intuitionistic fuzzy }(1,2)\text{-ideal of }S)=\nu_{A}(a).
\end{align*}
Similarly we can have $\mu_{B}(azxpa^{2})\geq\mu_{B}(a)$ and $\nu_{B}(azxpa^{2})\leq\nu_{B}(a).$ Then%
\begin{align*}
((\mu_{A})_{\beta\alpha}^{C}\circ(\mu_{B})_{\beta\alpha}^{C})(a)  &  =\underset
{a=mn}{\sup}\min\{(\mu_{A})_{\beta\alpha}^{C}(m),(\mu_{B})_{\beta\alpha}^{C}(n)\}\\
&  =\underset{a=mn}{\sup}\min\{\beta\cdot\mu_{A}(m)+\alpha,\beta\cdot
\mu_{B}(n)+\alpha\}\\
&  =\underset{a=mn}{\beta\cdot\sup}\min\{\mu_{A}(m),\mu_{B}(n)\}+\alpha\\
&  \geq\beta\cdot\min\{\mu_{A}(axypa^{2}),\mu_{B}(azxpa^{2})\}+\alpha\\
& (\text{since }a =(axypa^{2})(azxpa^{2}))\\
&  \geq\beta\cdot\min\{\mu_{A}(a),\mu_{B}(a)\}+\alpha(\text{since }A\text{ and
}B\text{ are intuitionistic}\\
& \text{fuzzy }(1,2)\text{-ideals} \text{ of }S)=\min\{\beta\cdot\mu_{A}(a)+\alpha,\beta\cdot\mu_{B}(a)+\alpha\}\\
&  =\min\{(\mu_{A})_{\beta\alpha}^{C}(a),\ (\mu_{B})_{\beta\alpha}^{C}(a)\}=((\mu_{A})_{\beta\alpha}^{C}\cap(\mu_{B})_{\beta\alpha}^{C})(a)
\end{align*} and%
\begin{align*}
((\nu_{A})_{\beta\alpha}^{C}\circ(\nu_{B})_{\beta\alpha}^{C})(a)  &  =\underset
{a=mn}{\inf}\max\{(\nu_{A})_{\beta\alpha}^{C}(m),(\nu_{B})_{\beta\alpha}^{C}(n)\}\\
&  =\underset{a=mn}{\inf}\max\{\beta\cdot\nu_{A}(m)-\alpha,\beta\cdot
\nu_{B}(n)-\alpha\}\\
&  =\underset{a=mn}{\beta\cdot\inf}\max\{\nu_{A}(m),\nu_{B}(n)\}-\alpha\\
&  \leq\beta\cdot\max\{\nu_{A}(axypa^{2}),\nu_{B}(azxpa^{2})\}-\alpha(\text{since }a =(axypa^{2})(azxpa^{2}))\\
&  \leq\beta\cdot\max\{\nu_{A}(a),\nu_{B}(a)\}-\alpha(\text{since }A\text{ and
}B\text{ are intuitionistic}\\
&\text{ fuzzy }(1,2)\text{-ideals of }S)=\max\{\beta\cdot\nu_{A}(a)-\alpha,\beta\cdot\nu_{B}(a)-\alpha\}\\
&  =\max\{(\nu_{A})_{\beta\alpha}^{C}(a),\ (\nu_{B})_{\beta\alpha}^{C}(a)\}=((\nu_{A})_{\beta\alpha}^{C}\cup(\nu_{B})_{\beta\alpha}^{C})(a).
\end{align*}

Hence $A_{\beta\alpha}^{C}$ $\circ B_{\beta\alpha}^{C}\supset A_{\beta\alpha}^{C}\cap B_{\beta\alpha}^{C}$ and $(1)$ follows. Similarly we can prove that $B_{\beta\alpha}^{C}$ $\circ A_{\beta\alpha}^{C}\supset A_{\beta\alpha}^{C}\cap B_{\beta\alpha}^{C}$. Combining these two results we can obtain $(2).$ Hence the proof.
\end{proof}

The following proposition can be proved by routine verification.

\begin{proposition}
Let $A=(\mu_{A},\nu_{A})$ be an intuitionistic fuzzy right ideal and $B=(\mu_{B},\nu_{B})$ be an intuitionistic fuzzy left ideal of a semigroup $S$. Then $A\circ B\subseteq A\cap B.$
\end{proposition}

\begin{proposition}
Let $S$ be a regular semigroup, $A=(\mu_{A},\nu_{A})$ be an intuitionistic fuzzy right ideal and $B=(\mu_{B},\nu_{B})$ be an intuitionistic fuzzy left ideal of $S$. Then $A\circ B\supseteq A\cap B.$
\end{proposition}

\begin{proof}
Let $c\in S.$ Then there exists an element $x\in S$ such that $c=cxc($since $S$ is regular$).$ Then%
\begin{align*}
(\mu_{A}\circ\mu_{B})(c)  &  =\underset{c=uv}{\sup}\{\min\{\mu
_{A}(u),\mu_{B}(v)\}\}\\
&  \geq\min\{\mu_{A}(c),\mu_{B}(c)\}=(\mu_{A}\cap\mu_{B})(c)
\end{align*}
and
\begin{align*}
(\nu_{A}\circ\nu_{B})(c)  &  =\underset{c=uv}{\inf}\{\max\{\nu
_{A}(u),\nu_{B}(v)\}\}\\
&  \leq\max\{\nu_{A}(c),\nu_{B}(c)\}=(\nu_{A}\cup\nu_{B})(c).
\end{align*}
Hence $A\circ B\supseteq A\cap B.$
\end{proof}

\begin{definition}
\cite{M} Let $S$ be a semigroup. Let $A$ and $B$ be subsets of $S.$ Then the {\it multiplication} of $A$ and $B$ is defined as \ $AB=\{ab\in S:a\in A$ and $b\in B\}.$
\end{definition}

\begin{theorem}
$\cite{M}$ A semigroup $S$ is regular if and only if $R\cap L=RL$ for every right ideal $R$ and every left ideal $L$ of $S.$
\end{theorem}

\begin{theorem}
For a semigroup $S$, the following conditions are equivalent:

$(1)$ $S$ is regular,

$(2)$ $A\circ B= A\cap B$ for any intuitionistic fuzzy right ideal $A=(\mu_{A},\nu_{A})$ and any intuitionistic fuzzy left ideal $B=(\mu_{B},\nu_{B})$ of $S$.
\end{theorem}

\begin{proof}
$(1)\Rightarrow(2):$ Let $S$ be a regular semigroup. Then by Proposition $3.27$ and $3.28,$ we have $A\circ B= A\cap B.$

$(2)\Rightarrow(1):$ Let $L$ and $R$ be respectively a left ideal and a right
ideal of $S$ and $x\in R\cap L.$ Then $x\in R$ and $x\in L.$ Hence $(\chi
_{L}(x),\overline{\chi}_{L}(x))=(\chi_{R}(x),\overline{\chi}_{R}(x))=(1,0)($where $\chi_{L}(x)$ and $\chi_{R}(x)$ are respectively the characteristic functions of $L$ and $R$$)$. Then%
\begin{align*}
(\chi_{R}\cap\chi_{L})(x)=\min\{\chi_{R}(x),\chi_{L}(x)\}=1\text{ and }(\overline{\chi}_{R}\cup\overline{\chi}_{L})(x)=\max\{\overline{\chi}_{R}(x),\overline{\chi}_{L}(x)\}=0.
\end{align*}
Now by Theorem $3.15,$ $(\chi_{L},\overline{\chi}_{L})$ and $(\chi_{R},\overline{\chi}_{R})$ are
respectively an intuitionistic fuzzy left ideal and an intuitionistic fuzzy right ideal of $S$. Hence by the hypothesis we have%
\begin{align*}
(\chi_{R}\circ\chi_{L})(x)=1,i.e.,\underset{x=yz}{\sup}[\min\{\chi_{R}(y),\chi_{L}(z)\}:y,z\in
S]=1
\end{align*}and%
\begin{align*}
(\overline{\chi}_{R}\circ\overline{\chi}_{L})(x)=0,i.e.,\underset{x=yz}{\inf}[\max\{\overline{\chi}_{R}(y),\overline{\chi}_{L}(z)\}:y,z\in
S]=0
\end{align*}
This implies that there exist some $r,s\in S$ such that $x=rs$ and
$\chi_{R}(r)=1=\chi_{L}(s)$ and $\overline{\chi}_{R}(r)=0=\overline{\chi}_{L}(s)$. Hence $r\in R$ and $s\in L.$ Hence $x\in RL.$
Thus $R\cap L\subseteq RL.$ Also $RL\subseteq R\cap L.$ Consequently,
$RL=R\cap L.$ Hence $S$ is regular.
\end{proof}

\begin{theorem}
For a semigroup $S$, the following conditions are equivalent:

$(1)$ $S$ is regular,

$(2)$ $A_{\beta\alpha}^{C}$ $\circ B_{\beta\alpha}^{C}=A_{\beta\alpha}^{C}\cap B_{\beta\alpha}^{C}$ for any intuitionistic fuzzy right ideal $A=(\mu_{A},\nu_{A})$ and any intuitionistic fuzzy left ideal $B=(\mu_{B},\nu_{B})$ of $S$.
\end{theorem}

\begin{proof}
$(1)\Rightarrow(2):$ Let $S$ be a regular semigroup, $A=(\mu_{A},\nu_{A})$ be an intuitionistic fuzzy right ideal
and $B=(\mu_{B},\nu_{B})$ be an intuitionistic fuzzy left ideal of $S$. Then by Theorem $3.9,$ $A_{\beta\alpha}^{C}$ is an intuitionistic fuzzy right ideal and $B_{\beta\alpha}^{C}$ is an intuitionistic fuzzy left ideal of $S.$ Hence by Theorem $3.31,$ $A_{\beta\alpha}^{C}$ $\circ B_{\beta\alpha}^{C}=A_{\beta\alpha}^{C}\cap B_{\beta\alpha}^{C}.$

$(2)\Rightarrow(1):$ Let $L$ and $R$ be respectively a left ideal and a right ideal of $S$ and $x\in R\cap L.$ Then $x\in R$ and $x\in L.$ Hence $(\chi_{L}(x),\overline{\chi}_{L}(x))=(\chi_{R}(x),\overline{\chi}_{R}(x))=(1,0)($where $\chi_{L}(x)$ and $\chi_{R}(x)$ are respectively the characteristic functions of $L$ and $R$$)$. Since $x\in L,$ $\overline{\chi}_{L}(x)=0.$ This implies that $\beta\cdot\overline{\chi}_{L}(x)=0.$ Hence $\inf\{\beta\cdot\overline{\chi}_{L}(y):y\in S\}=0.$ Since $\alpha\in[0,\inf\{\beta\cdot\overline{\chi}_{L}(y):y\in S\}],$ $\alpha=0.$ Similarly, we can show that if $x\in R,$ then also $\alpha=0.$ Then $((\chi_{L})_{\beta\alpha}^{C}(x),(\overline{\chi}_{L})_{\beta\alpha}^{C}(x))=((\chi_{R})_{\beta\alpha}^{C}(x),(\overline{\chi}_{R})_{\beta\alpha}^{C}(x))=(\beta, 0).$ Thus%
\begin{align*}
((\chi_{R})_{\beta\alpha}^{C}\cap(\chi_{L})_{\beta\alpha}^{C})(x)  &
=\min\{(\chi_{R})_{\beta\alpha}^{C}(x),(\chi_{L})_{\beta\alpha}^{C}(x)\}\\
&  =\beta.
\end{align*} and%
\begin{align*}
((\overline{\chi}_{R})_{\beta\alpha}^{C}\cup(\overline{\chi}_{L})_{\beta\alpha}^{C})(x)  &
=\max\{(\overline{\chi}_{R})_{\beta\alpha}^{C}(x),(\overline{\chi}_{L})_{\beta\alpha}^{C}(x)\}\\
&  =0.
\end{align*}
Now by Theorem $3.15$, $(\chi_{L},\overline{\chi}_{L})$ and $(\chi_{R},\overline{\chi}_{R})$ are
respectively an intuitionistic fuzzy left ideal and an intuitionistic fuzzy right ideal of $S$. Hence by
Theorem $3.9,$ $((\chi_{L})_{\beta\alpha}^{C},(\overline{\chi}_{L})_{\beta\alpha}^{C})$ and $((\chi_{R})_{\beta\alpha
}^{C},(\overline{\chi}_{R})_{\beta\alpha}^{C})$ are respectively intuitionistic fuzzy left ideal and intuitionistic fuzzy right ideal of $S.$ This
together with the hypothesis gives%
\begin{align*}
((\chi_{R})_{\beta\alpha}^{C}\circ(\chi_{L})_{\beta\alpha}^{C})(x)  &
=\beta\\
i.e.,\underset{x=yz}{\sup}[\min\{(\chi_{R})_{\beta\alpha}^{C}(y),(\chi
_{L})_{\beta\alpha}^{C}(z)\}  &  :y,z\in S]=\beta\\
i.e.,\underset{x=yz}{\sup}[\min\{\beta\cdot\chi_{R}(y)+\alpha,\beta\cdot
\chi_{L}(z)+\alpha\}  &  :y,z\in S]=\beta\\
i.e.,\beta.\underset{x=yz}{\sup}[\min\{\chi_{R}(y),\chi_{L}(z)\}  &  :y,z\in
S]+\alpha=\beta\\
i.e.,\beta.\underset{x=yz}{\sup}[\min\{\chi_{R}(y),\chi_{L}(z)\}  &  :y,z\in
S]=\beta(\text{since }\alpha=0)
\end{align*} and%
\begin{align*}
((\overline{\chi}_{R})_{\beta\alpha}^{C}\circ(\overline{\chi}_{L})_{\beta\alpha}^{C})(x)  &
=0\\
i.e.,\underset{x=yz}{\inf}[\max\{(\overline{\chi}_{R})_{\beta\alpha}^{C}(y),(\overline{\chi}
_{L})_{\beta\alpha}^{C}(z)\}  &  :y,z\in S]=0\\
i.e.,\underset{x=yz}{\inf}[\max\{\beta\cdot\overline{\chi}_{R}(y)-\alpha,\beta\cdot
\overline{\chi}_{L}(z)-\alpha\}  &  :y,z\in S]=0\\
i.e.,\beta.\underset{x=yz}{\inf}[\max\{\overline{\chi}_{R}(y),\overline{\chi}_{L}(z)\}  &  :y,z\in
S]-\alpha=0
\end{align*}
\begin{align*}
i.e.,\beta.\underset{x=yz}{\inf}[\max\{\overline{\chi}_{R}(y),\overline{\chi}_{L}(z)\}  &  :y,z\in
S]=0(\text{since }\alpha=0)
\end{align*}
Hence $\underset{x=yz}{\sup}[\min\{\chi_{R}(y),$ $\ \chi_{L}(z)\}:y,z\in
S]=1$ and $\underset{x=yz}{\inf}[\max\{\overline{\chi}_{R}(y),$ $\ \overline{\chi}_{L}(z)\}:y,z\in
S]=0$ This implies that there exist some $r,s\in S$ such that $x=rs$ and
$\chi_{R}(r)=1=\chi_{L}(s),$ $\overline{\chi}_{R}(r)=0=\overline{\chi}_{L}(s)$. Hence $r\in R$ and $s\in L$ whence $x\in RL.$
Thus $R\cap L\subseteq RL.$ Also $RL\subseteq R\cap L.$ Consequently,
$RL=R\cap L.$ Hence $S$ is regular.
\end{proof}

\begin{remark}
If we put $\beta=1($respectively $\alpha=0)$ in intuitionistic fuzzy magnified translation
then it reduces to intuitionistic fuzzy translation$($respectively intuitionistic fuzzy multiplication$).$
Consequently analogues of Theorems $3.4$-$3.6, 3.8$-$3.11,$ Theorem $3.14$-$3.15,$ Theorem $3.17,3.19$-$3.20, 3.22, 3.25$-$3.26, 3.32$ and Corollary $3.12$ follow easily in intuitionistic fuzzy translation and intuitionistic fuzzy multiplication.
\end{remark}


\end{document}